\documentclass[preprint,12pt,a4paper,openright]{elsarticle}

\usepackage{setspace} 

\makeatletter
\def\ps@pprintTitle{%
	\let\@oddhead\@empty
	\let\@evenhead\@empty
	\def\@oddfoot{\reset@font\hfil\thepage\hfil}
	\let\@evenfoot\@oddfoot
}
\makeatother
\usepackage{upgreek}
\usepackage{fixmath}
\usepackage{mathptmx}	
\usepackage[utf8]{inputenc}
\usepackage[english]{babel}
\usepackage{amsthm,amssymb,amsmath}
\usepackage{bigints}
\usepackage{multirow}
\usepackage{xcolor}
\usepackage{mathrsfs}
\usepackage{graphicx}
\usepackage{epsfig}
\usepackage{caption}
\usepackage{subcaption}
\usepackage{lscape} 
\usepackage{enumitem}

\usepackage{rotating}

\usepackage{chngcntr} 
\counterwithin{table}{section} 
\counterwithin{figure}{section} 

\usepackage{natbib}
\setcitestyle{authoryear,open={(},close={)}} 

\usepackage[colorlinks=true]{hyperref}
\hypersetup{
	colorlinks=true,
	linkcolor=blue,  
	citecolor=blue,
}

\newtheorem{eg}{Example}[section]
\newtheorem{ceg}{Counter Example}[section]
\newtheorem{thm}{Theorem}[section]

\numberwithin{equation}{section}

\title{Further Results on the Bivariate Semi-parametric Singular Family of Distributions}
\author{{Durga Vasudevan and G. Asha\footnote{Corresponding author email: asha@cusat.ac.in}}\\
	{\textit{Department of Statistics}}\\
	{\textit{Cochin University of Science and Technology, Cochin, Kerala, India-682022}}\\
	{Email: \textit{durgavvn10@gmail.com, asha@cusat.ac.in}}}

\begin{document}
	
	\begin{frontmatter}
		\begin{abstract}
			General classes of bivariate distributions are well studied in literature. Most of these classes are proposed via a copula formulation or extensions of some characterisation properties in the univariate case. In \citet{kundu2022bivariate} we see one such semi-parametric family useful to model bivariate data with ties. This model is a general semi-parametric model with a baseline. In this paper we present a characterisation property of this class of distributions in terms of a functional equation. The general solution to this equation is explored. Necessary and sufficient conditions under which the solution becomes a bivariate distribution is investigated.
			
		\end{abstract}
	\end{frontmatter}
	
	\section{Introduction}
	Analysing bivariate datasets are very challenging as it involves association between two variables. Data with ties are also common. The Marshall-Olkin copula defined by
	\begin{equation}\label{mo_copula}
		C(u_1,u_2)=\begin{cases}
			u_1^{1-a} u_2~;~u_1^a \geq u_2^b\\
			u_1 u_2^{1-b}~;~u_1^a \leq u_2^b
		\end{cases}
	\end{equation}
	where $0 < a,b < 1$ and $0 < u_i < 1;~i=1,2$ is apt in modelling association of these types (see \citet{nelsen2007introduction}). Many authors have proposed very general classes of bivariate distributions and studied their properties. In \eqref{mo_copula}, when $u_i$ is exponential($\theta_i+\theta_3$) for $i=1,2$ and  $a=\frac{\theta_3}{\theta_1 + \theta_3}$ and $b=\frac{\theta_3}{\theta_2 + \theta_3}$, the well studied Marshall-Olkin bivariate exponential distribution with a singular component is obtained as
	\begin{equation}\label{mo_distribution}
		\bar{F}(x_1,x_2)=
		\begin{cases}
			e^{-(\theta_1 + \theta_3) x_1 - \theta_2 x_2} &~;~x_{1} \geq x_{2}\\
			e^{-\theta_1 x_1 - (\theta_2 + \theta_3) x_2} &~;~x_{1} \leq x_{2},
		\end{cases}
	\end{equation}
	where $\theta_i > 0;~i=1,2,3$.
	Here $(X_1,X_2)$ is a bivariate random vector with $\bar{F}(x_1,x_2)=P(X_1 > x_1, X_2 > x_2)$ with support $I_{(0, \infty) \times (0, \infty)}=\{ (x_1,x_2);~0 < x_1 < \infty,~0 < x_2 < \infty\}$. The distribution in \eqref{mo_distribution} is characterised by the functional equation popularly known as bivariate lack of memory property (BLMP) given by
	\begin{equation*}
		\bar{F}(x_1 + t, x_2 + t)=\bar{F}(x_1,x_2) \bar{F}(t,t)
	\end{equation*}
	for all $x_1,x_2,t \geq 0$. \citet{kolev2018weak} worked on the Marshall-Olkin's bivariate exponential distribution thus providing a weak version of the BLMP which can be used to construct bivariate distributions having a singularity component along arbitrary line through the origin. \citet{lin2019bivariate} studied the moment generating function, product moments and dependence structure of the bivariate distributions satisfying BLMP. Many authors have studied the BLMP under different perspectives (\citet{pinto2015sibuya} and \citet{kolev2018functional}). In fact the BLMP has been translated as a property of arbitrary bivariate continuous distributions (\citet{galambos2006characterizations}).
	\par If $u_i$ is Weibull($\theta_i + \theta_3,\alpha$) for $i=1,2$, the Marshall-Olkin copula in \eqref{mo_copula} gives a bivariate distribution with a singular component whose survival function is
	\begin{equation}\label{mo_weibull}
		\bar{F}(x_1,x_2)=
		\begin{cases}
			e^{-(\theta_1 + \theta_3) x_1^{\alpha} - \theta_2 x_2^{\alpha}} &~;~x_{1} \geq x_{2}\\
			e^{-\theta_1 x_1^{\alpha} - (\theta_2 + \theta_3) x_2^{\alpha}} &~;~x_{1} \leq x_{2},
		\end{cases}
	\end{equation}
	where $\alpha,~\theta_i > 0;~i=1,2,3$, characterised by the functional equation,
	\begin{equation*}
		\bar{F}((x_1^{\alpha} + t^{\alpha})^{\frac{1}{\alpha}}, (x_2^{\alpha} + t^{\alpha})^{\frac{1}{\alpha}})=\bar{F}(x_1,x_2) \bar{F}(t,t),
	\end{equation*}
	for all $x_1,x_2,t \geq 0$.
	\par If $u_i$ is Pareto($\theta_i + \theta_3$) for $i=1,2$, the Marshall-Olkin copula in \eqref{mo_copula} gives a bivariate distribution with a singular component whose survival function is
	\begin{equation}\label{mo_pareto}
		\bar{F}(x_1,x_2)=
		\begin{cases}
			(\frac{1}{x_1})^{\theta_1 + \theta_3} (\frac{1}{x_2})^{\theta_2} &~;~x_{1} \geq x_{2}\\
			(\frac{1}{x_1})^{\theta_1} (\frac{1}{x_2})^{\theta_2 + \theta_3} &~;~x_{1} \leq x_{2},
		\end{cases}
	\end{equation}
	where $\theta_i > 0;~i=1,2,3$, characterised by the functional equation,
	\begin{equation*}
		\bar{F}(x_1 t , x_2 t)=\bar{F}(x_1,x_2) \bar{F}(t,t),
	\end{equation*}
	for all $x_1,x_2,t \geq 0$.
	\par If the marginals belong to the proportional hazard (PH) class given by $u_i=\bar{F}_i(x_i)=[\bar{F}_0(x_i)]^{\theta_i + \theta_3}$ for $i=1,2$, then \eqref{mo_copula} becomes the copula associated with the bivariate proportional hazard class (BPHC) mentioned in \citet{kundu2022bivariate}
	whose bivariate survival function is
	\begin{equation}\label{kundu_surv1}
		\bar{F}(x_1,x_2)=
		\begin{cases}
			[\bar{F}_0(x_1)]^{\theta_1 + \theta_3} [\bar{F}_0(x_2)]^{\theta_2} &~;~x_{1} \geq x_{2}\\
			[\bar{F}_0(x_1)]^{\theta_1} [\bar{F}_0(x_2)]^{\theta_2 + \theta_3} &~;~x_{1} \leq x_{2},
		\end{cases}
	\end{equation}
	where $\bar{F}_0(x)=P(X > x)$ is the baseline distribution with $\bar{F}_0(x)=0$ for $x \leq 0$ without loss of generality and $\theta_i > 0;~i=1,2,3$. In this work, we study the functional equation characterised by singular family of bivariate distributions in \eqref{kundu_surv1}.
	The rest of the paper is organised as below.
	\par In Section \ref{fun_eqn}, we propose a class of distributions as a general solution to the functional equation characterising \eqref{kundu_surv1}. This functional equation characterises a very general class of distributions which includes \eqref{mo_distribution}, \eqref{mo_weibull}, \eqref{mo_pareto} and \eqref{kundu_surv1}. It is observed through a counter example that the general solutions of this functional equation need not necessarily generate a bivariate distribution function. But conditions on the marginals ensures that the solution is a bivariate distribution. In Section \ref{sec_ns_conditions}, necessary and sufficient conditions are derived for a univariate distribution to be a marginal for the bivariate distribution belonging to this class. In Section \ref{sec_construction}, we translate these in terms of failure rates for the ease of constructing bivariate distributions satisfying the functional equation.
	
	\section{Functional equation}\label{fun_eqn}
	In this section, we propose a characterisation of \eqref{kundu_surv1} in terms of a functional equation, the general solution of which includes many bivariate singular distributions including \eqref{kundu_surv1}.
	
	\begin{thm}
		For some baseline survival function $\bar{F}_0(\cdot)$, the functional equation
		\begin{equation}\label{kundu_fe_gen}
			\bar{F}(\bar{F}_0^{-1}(\bar{F}_0(t)\bar{F}_0(x_1)), \bar{F}_0^{-1}(\bar{F}_0(t)\bar{F}_0(x_2)))=\bar{F}(x_1,x_2) \bar{F}(t,t)
		\end{equation}
		is satisfied if and only if \begin{equation}\label{kundu_surv_gen}
			\bar{F}(x_1,x_2)=
			\begin{cases}
				\bar{F_1} \Big( \bar{F}_0^{-1} \big( \frac{\bar{F}_0(x_1)}{\bar{F}_0(x_2)} \big) \Big) [\bar{F}_0(x_2)]^{\theta} &~;~x_{1} \geq x_{2}\\
				\bar{F_2} \Big( \bar{F}_0^{-1} \big( \frac{\bar{F}_0(x_2)}{\bar{F}_0(x_1)} \big) \Big) [\bar{F}_0(x_1)]^{\theta} &~;~x_{1} \leq x_{2},
			\end{cases}
		\end{equation}
		for some $\theta >0$ and $x_1,x_2 \geq 0$ where $\bar{F_i}(x_i)=\bar{F}(x_1,x_2)_{\vert x_{3-i}=0};~i=1,2$.
	\end{thm}
	
	\begin{proof}
		In \eqref{kundu_fe_gen}, for $x_1=x_2=x$, we get	\begin{equation}\label{kundu_fe1}
			\bar{F}(\bar{F}_0^{-1}(\bar{F}_0(t)\bar{F}_0(x)), \bar{F}_0^{-1}(\bar{F}_0(t)\bar{F}_0(x)))=\bar{F}(x,x) \bar{F}(t,t).
		\end{equation}
		Writing $\bar{F}(x,x)=H(x)$, \eqref{kundu_fe1} becomes,
		\begin{equation}\label{kundu_fe1_univ}
			H(\bar{F}_0^{-1}(\bar{F}_0(t)\bar{F}_0(x)))=H(x)H(t)
		\end{equation}
		which is of the form
		\begin{equation*}
			T(us)=T(s)T(u),
		\end{equation*}
		where $s=\bar{F}_0(x)$, $u=\bar{F}_0(t)$ and $H \circ \bar{F}_0^{-1}=T$. Hence from  (\citet{aczel1966lectures}, pg. 38), the most general solution of \eqref{kundu_fe1_univ} is
		$$H(t)=\bar{F}(t,t)=[\bar{F}_0(t)]^{\theta};~\theta >0.$$
		Therefore, \eqref{kundu_fe_gen} can be written as,
		\begin{equation*}\label{kundu_fe2}
			\bar{F}(\bar{F}_0^{-1}(\bar{F}_0(t)\bar{F}_0(x_1)), \bar{F}_0^{-1}(\bar{F}_0(t)\bar{F}_0(x_2)))=\bar{F}(x_1,x_2) [\bar{F}_0(t)]^{\theta}.
		\end{equation*}
		For $x_1 \geq x_2$,
		\begin{align*}
			\bar{F}(x_1,x_2)&=\bar{F}\Big( \bar{F}_0^{-1} \big( \bar{F}_0(x_2) \frac{\bar{F}_0(x_1)}{\bar{F}_0(x_2)} \big), \bar{F}_0^{-1}(\bar{F}_0(x_2)) \big)\\
			&=\bar{F}_1\Big( \bar{F}_0^{-1} \big( \frac{\bar{F}_0(x_1)}{\bar{F}_0(x_2)}\big) \Big) [\bar{F}_0(x_2)]^{\theta}.
		\end{align*}
		Arguing similarly for $x_1 \leq x_2$, \eqref{kundu_surv_gen} can be retrieved. The converse is direct.
	\end{proof}
	
	The next question of interest is whether this general equation $\bar{F}(x_1,x_2)$ in \eqref{kundu_surv_gen} represent proper bivariate distribution function? The following counter example shows it is not. It is further seen that arbitrary marginals do not admit a proper bivariate distribution.
	
	\begin{ceg}\label{kundu_counter_eg}
		\normalfont LFR-exponential: For the linear failure rate (LFR) marginals, $\bar{F_i}(x_i)=e^{-(x_i + \alpha x_i^2)};~x_i > 0, \alpha >0,~i=1,2$ and exponential baseline, $\bar{F}_0(x)=e^{-x};~x>0$, equation \eqref{kundu_fe_gen} reduces to
		\begin{equation}\label{kundu_fe_LFR_exp}
			\bar{F}(x_1 + t,x_2 + t)= \bar{F}(x_1,x_2) \bar{F}(t,t)
		\end{equation}
		with the corresponding solution \eqref{kundu_surv_gen} as
		\begin{equation}\label{kundu_counter_eg_surv}
			\bar{F}(x_1,x_2)=
			\begin{cases}
				e^{-(x_1-x_2)-\alpha (x_1-x_2)^2-\theta x_2} &~;~x_{1} \geq x_{2}\\
				e^{-(x_2-x_1)-\alpha (x_2-x_1)^2-\theta x_1} &~;~x_{1} \leq x_{2}.
			\end{cases}
		\end{equation}
		On closely examining $\bar{F}(x_1,x_2)$ we observe that $P(1 < X_1 < 2, 3 < X_2 < 5) \ngeq 0$ for the choice of $\alpha=1.5$ and $\theta=3$ disqualifying it as a bivariate probability survival function. The functional equation \eqref{kundu_fe_LFR_exp} has been well studied under this perspective in \citet{kulkarni2006characterizations}. Hence not all solution with arbitrary $\bar{F_i}(\cdot)$ and $\bar{F}_0(\cdot)$ give bivariate probability distribution. Imposing some restrictions on the marginals ensures this. Motivated by this we develop the conditions to be satisfied by the marginals so that $\bar{F}(x_1,x_2)$ is a bivariate survival function.
	\end{ceg}
	
	\section{Necessary and sufficient conditions for generating bivariate distributions}\label{sec_ns_conditions}
	In this section, we discuss the conditions to be satisfied by the univariate distributions for them to qualify as a marginal for the bivariate survival function as in \eqref{kundu_surv_gen}. Let $\Theta = (\theta, \mathbold{\alpha}, \mathbold{\beta})$, where $\mathbold{\alpha} = (\alpha_1,\alpha_2,\dots,\alpha_k)$ is the vector of parameters involved in baseline $\bar{F}_0(\cdot)$ and $\mathbold{\beta} = (\beta_1,\beta_2,\dots,\beta_m)$ is the vector of parameters involved in the marginal $\bar{F_i}(\cdot)$.
	
	\begin{thm}\label{kundu_mo_theorem}
		Let $F_i(x)$ be a distribution function with absolutely continuous density $f_i(x)$ for which $\lim_{x \rightarrow \infty} f_i(x)=0;~i=1,2$. The necessary and sufficient conditions for $\bar{F}(x_1,x_2)$ in \eqref{kundu_surv_gen} be a bivariate distribution is that
		\begin{enumerate}[label=(\roman*)]
			\item $\theta \leq u_1(\Theta) + u_2(\Theta) \leq 2 \theta$ \label{kundu_marg_densty1}
			\item $\frac{\partial}{\partial x_{3-i}} \ln \Big( -\frac{\partial}{\partial x_{i}} \bar{F_i} \Big( \bar{F}_0^{-1} \big( \frac{\bar{F}_0(x_i)}{\bar{F}_0(x_{3-i})} \big) \Big) \Big) \leq \theta r_0(x_{3-i});~i=1,2$, \label{kundu_marg_densty2}
		\end{enumerate} 
		where $u_i(\Theta)=\Big[f_i\Big( \bar{F}_0^{-1} \big( \frac{\bar{F}_0(x_i)}{\bar{F}_0(x_{3-i})} \big) \Big) \Big \vert \frac{\partial \bar{F}_0^{-1} \big( \frac{\bar{F}_0(x_i)}{\bar{F}_0(x_{3-i})} \big) }{\partial (- \ln \bar{F}_0(x_{3-i}))} \Big \vert \Big] _{\vert x_i=x_{3-i}};~i=1,2$.
	\end{thm}
	
	\begin{proof}
		$\bar{F}(x_1,x_2)$ will be a bivariate survival function if and only if it can be written as a convex mixture of $\bar{F}_{a}(x_1,x_2)$ and $\bar{F}_{s}(x_1,x_2)$ where $\bar{F}_{a}(x_1,x_2)$ is the absolutely continuous part and $\bar{F}_{s}(x_1,x_2)$ is the singular part and $\bar{F}_{a}(x_1,x_2)$ and $\bar{F}_{s}(x_1,x_2)$ are survival functions.
		Observe, for $0 \leq \alpha \leq 1$,
		\begin{align}
			\nonumber \frac{\partial^2 \bar{F}(x_1,x_2)}{\partial x_1 \partial x_2}&= \alpha f_a(x_1,x_2) \\ \label{kundu_absdensty}&=
			\begin{cases}
				[\bar{F}_0(x_2)]^{\theta} \big[ \frac{\partial^2 \bar{F}_{1} \big( \bar{F}_0^{-1} \big( \frac{\bar{F}_0(x_1)}{\bar{F}_0(x_2)} \big) \big) }{\partial x_1 \partial x_2} - \theta \frac{f_0(x_2)}{\bar{F}_0(x_2)} \frac{\partial \bar{F}_1 \big( \bar{F}_0^{-1} \big( \frac{\bar{F}_0(x_1)}{\bar{F}_0(x_2)} \big) \big)}{\partial x_1} \big] &~;~x_{1} \geq x_{2}\\
				[\bar{F}_0(x_1)]^{\theta} \big[ \frac{\partial^2 \bar{F}_{2} \big( \bar{F}_0^{-1} \big( \frac{\bar{F}_0(x_2)}{\bar{F}_0(x_1)} \big) \big) }{\partial x_1 \partial x_2} - \theta \frac{f_0(x_1)}{\bar{F}_0(x_1)} \frac{\partial \bar{F}_2 \big( \bar{F}_0^{-1} \big( \frac{\bar{F}_0(x_2)}{\bar{F}_0(x_1)} \big) \big)}{\partial x_2} \big] &~;~x_{1} \leq x_{2}.
			\end{cases}
		\end{align}
		Also,
		\begin{align*}
			\int_{x_1 \geq x_2} \alpha f_a(x_1,x_2)
			=&1-\frac{1}{\theta} \bigg[f_1 \Big( \bar{F}_0^{-1} \big( \frac{\bar{F}_0(x_1)}{\bar{F}_0(x_{2})} \big) \Big) \bigg \vert \frac{\partial \bar{F}_0^{-1} \big( \frac{\bar{F}_0(x_1)}{\bar{F}_0(x_2)} \big) }{\partial (- \ln \bar{F}_0(x_2))} \bigg \vert \bigg] _{\vert x_1=x_2}
		\end{align*}
		and
		$$\int_{x_1 \leq x_2} \alpha f_a(x_1,x_2) = 1-\frac{1}{\theta} \bigg[f_2 \Big( \bar{F}_0^{-1} \big( \frac{\bar{F}_0(x_2)}{\bar{F}_0(x_{1})} \big) \Big) \bigg \vert \frac{\partial \bar{F}_0^{-1} \big( \frac{\bar{F}_0(x_1)}{\bar{F}_0(x_2)} \big) }{\partial (- \ln \bar{F}_0(x_2))} \bigg \vert \bigg] _{\vert x_2=x_1}.$$
		Writing, $$u_i(\Theta)=\Big[f_i\Big( \bar{F}_0^{-1} \big( \frac{\bar{F}_0(x_i)}{\bar{F}_0(x_{3-i})} \big) \Big) \Big \vert \frac{\partial \bar{F}_0^{-1} \big( \frac{\bar{F}_0(x_i)}{\bar{F}_0(x_{3-i})} \big) }{\partial (- \ln \bar{F}_0(x_{3-i}))} \Big \vert \Big] _{\vert x_i=x_{3-i}};~i=1,2$$ we have,
		\begin{eqnarray*}
			\int_{x_1 \geq x_2} \alpha f_a(x_1,x_2) &= 1-\frac{1}{\theta}u_1(\Theta) \\
			\text{and}~\int_{x_1 \leq x_2} \alpha f_a(x_1,x_2) &= 1-\frac{1}{\theta}u_2(\Theta)
		\end{eqnarray*}
		so that
		\begin{equation}\label{kundu_absdensty_alpha}
			\alpha=2- \frac{1}{\theta} (u_1(\Theta) + u_2(\Theta)).
		\end{equation}
		Hence the absolutely continuous part has density $f_a(x_1,x_2)$ given by \eqref{kundu_absdensty} and \eqref{kundu_absdensty_alpha}. Also, we have
		$$\bar{F}_a(x,x)=[\bar{F}_0(x)]^{\theta}.$$
		But $\bar{F}(x,x)=[\bar{F}_0(x)]^{\theta}$, so that
		$\bar{F}_s(x,x)=\frac{\bar{F}(x,x)-\alpha \bar{F}_a(x,x)}{1-\alpha}=[\bar{F}_0(x)]^{\theta}$. Therefore, $\bar{F}$ is a valid survival function if
		\begin{enumerate}[label=(\roman*)]
			\item $\bar{F}$ is a convex mixture of $F_a$ and $F_s$. From \eqref{kundu_absdensty_alpha}, since $0 \leq \alpha \leq 1$, we get $$\theta \leq u_1(\Theta) + u_2(\Theta) \leq 2 \theta.$$
			\item $\bar{F}_a$ is a valid survival function. From \eqref{kundu_absdensty}, since $f_a(x_1,x_2) \geq 0$, we get for $x_i > x_{3-i};~i=1,2$,
			$$\frac{\partial}{\partial x_{3-i}} \ln \Big( -\frac{\partial}{\partial x_{i}} \bar{F_i} \Big( \bar{F}_0^{-1} \big( \frac{\bar{F}_0(x_i)}{\bar{F}_0(x_{3-i})} \big) \Big) \Big) \leq \theta r_0(x_{3-i}).$$ 
		\end{enumerate} 
		This completes the proof.
	\end{proof}
	
	\textbf{Counter Example \ref{kundu_counter_eg} (contd.)} LFR-exponential: Consider again the bivariate function $\bar{F}(x_1,x_2)$ given in \eqref{kundu_counter_eg_surv}. Condition \ref{kundu_marg_densty2} of Theorem \ref{kundu_mo_theorem} becomes
	\begin{equation*}
		2 \alpha (x_i-x_{3-i}) - \frac{1}{x_i - x_{3-i}} + 1 \leq \theta;~x_i > x_{3-i},~i=1,2,
	\end{equation*}
	which is not satisfied for the choice of $\alpha=1.5,~\theta=3,~x_1=5$ and $x_2=3$. This again shows that a bivariate distribution cannot be formed with marginal being the linear failure rate model and baseline as exponential distribution in \eqref{kundu_surv_gen}. However any proportional hazard distribution qualify as marginals in \eqref{kundu_surv_gen}. This is discussed in detail in the next theorem.

	\begin{thm}
		The survival function $\bar{F}(x_1,x_2)$ with the marginals belonging to the PH class satisfies the functional equation \eqref{kundu_fe_gen} if and only if it is of the form
		\begin{align}\label{kundu_surv_gen2}
			\bar{F}(x_1,x_2)&=
			\begin{cases}
				[\bar{F}_0(x_1)]^{\delta_1} [\bar{F}_0(x_2)]^{\theta - \delta_1} &~;~x_{1} \geq x_{2}\\
				[\bar{F}_0(x_1)]^{\theta - \delta_2} [\bar{F}_0(x_2)]^{\delta_2} &~;~x_{1} \leq x_{2}.
			\end{cases}
		\end{align}
		for some $\theta,~\delta_1,~\delta_2 >0$ with $\delta_i < \theta;~i=1,2$.
	\end{thm}
	
	\begin{proof}
		Suppose \eqref{kundu_fe_gen} is satisfied. Then, the general solution of $\bar{F}(x_1,x_2)$ is given by \eqref{kundu_surv_gen}. Let the marginals belong to the PH class so that $\bar{F_i}(x_i)=[\bar{F}_0(x_i)]^{\delta_i};~\delta_i>0,~i=1,2$. Then, the bivariate function in \eqref{kundu_surv_gen} reduces to \eqref{kundu_surv_gen2}. Now, it remains to show that $\bar{F}(x_1,x_2)$ is a valid bivariate distribution. Proceeding as in Theorem \ref{kundu_mo_theorem},
		\begin{align}
			\nonumber \frac{\partial^2 \bar{F}(x_1,x_2)}{\partial x_1 \partial x_2}&= \alpha f_a(x_1,x_2) \\ \label{kundu_absdensty1}&=
			\begin{cases}
				\delta_1 (\theta - \delta_1) f_0(x_1) f_0(x_2) (\bar{F}_0(x_1))^{\delta_1 - 1} (\bar{F}_0(x_2))^{\theta - \delta_1 - 1} &~;~x_{1} \geq x_{2}\\
				\delta_2 (\theta - \delta_2) f_0(x_1) f_0(x_2) (\bar{F}_0(x_1))^{\theta - \delta_2 - 1} (\bar{F}_0(x_2))^{\delta_2 - 1} &~;~x_{1} \leq x_{2}.
			\end{cases}
		\end{align}
		Then,
		\begin{align*}
			\int_{x_1 \geq x_2} \alpha f_a(x_1,x_2) 
			=&\frac{\theta - \delta_1}{\theta} \\
			\text{and}~ \int_{x_1 \leq x_2} \alpha f_a(x_1,x_2) 
			=&\frac{\theta - \delta_2}{\theta}
		\end{align*}
		so that $\alpha=\frac{2 \theta - \delta_1 - \delta_2}{\theta}$.
		Hence, the absolute continuous part has density given by,
		\begin{equation}\label{kundu_abs_densty1}
			f_a(x_1,x_2)= \frac{\theta}{2 \theta - \delta_1 - \delta_2} \times
			\begin{cases}
				\delta_1 (\theta - \delta_1) f_0(x_1) f_0(x_2) (\bar{F}_0(x_1))^{\delta_1 - 1} (\bar{F}_0(x_2))^{\theta - \delta_1 - 1} &~;~x_{1} \geq x_{2}\\
				\delta_2 (\theta - \delta_2) f_0(x_1) f_0(x_2) (\bar{F}_0(x_1))^{\theta - \delta_2 - 1} (\bar{F}_0(x_2))^{\delta_2 - 1} &~;~x_{1} \leq x_{2}.
			\end{cases}
		\end{equation}
		We have, $\bar{F}_a(x_1,x_2)=[\bar{F}_0(x)]^{\theta}$. But, $\bar{F}(x,x)=[\bar{F}_0(x)]^{\theta}$ so that $\bar{F}_s(x_1,x_2)=\frac{\bar{F}(x,x)-\alpha \bar{F}_a(x,x)}{1-\alpha}=[\bar{F}_0(x)]^{\theta}$. From \ref{kundu_marg_densty1} and \ref{kundu_marg_densty2} as in Theorem \ref{kundu_mo_theorem}, we have $\theta \leq \delta_1 + \delta_2 \leq 2 \theta$ which is true since $\delta_i < \theta;~i=1,2$. Also, $\bar{F}_a$ is a valid survival function implying $f_a(x_1,x_2) \geq 0$ which is true for any choice of the baseline distribution $\bar{F}_0(\cdot)$. The converse is direct.
	\end{proof}
	
	
	\section{Construction of bivariate distributions satisfying functional equation \eqref{kundu_fe_gen}}\label{sec_construction}
	The functional equation \eqref{kundu_fe_gen} can be equivalently represented in terms of the hazard gradient. The discussions above translated in terms of the hazard gradient would make it computationally easier to check if marginals qualify for the function \eqref{kundu_surv_gen} to be a bivariate survival function. Accordingly in this section we visit the functional equation in terms of its hazard gradient. The hazard gradient vector $r(x_1,x_2)$ (\citet{johnson1975vector}) is defined as
	\begin{align*}
		r(x_1,x_2)&=(r_1(x_1,x_2),r_2(x_1,x_2))\\
		&=\Big( -\frac{\partial \ln \bar{F}(x_1,x_2)}{\partial x_1},-\frac{\partial \ln \bar{F}(x_1,x_2)}{\partial x_2} \Big)\\
		&=\Big( \frac{\partial}{\partial x_1} R(x_1,x_2),  \frac{\partial}{\partial x_2} R(x_1,x_2) \Big),
	\end{align*}
	where $R(x_1,x_2)=- \ln \bar{F}(x_1,x_2)$.
	
	\par Note that the functional equation \eqref{kundu_fe_gen} can be equivalently expressed as
	\begin{equation*}
		R(\bar{F}_0^{-1}(\bar{F}_0(t)\bar{F}_0(x_1)), \bar{F}_0^{-1}(\bar{F}_0(t)\bar{F}_0(x_2)))=R(x_1,x_2)-\theta \ln \bar{F}_0(t).
	\end{equation*}
	Differentiating with respect to $t$, we get,
	\begin{align*}
		&\sum_{i=1}^{2} \frac{\partial \big [R(\bar{F}_0^{-1}(\bar{F}_0(t)\bar{F}_0(x_1)), \bar{F}_0^{-1}(\bar{F}_0(t)\bar{F}_0(x_2)))]}{\partial (\bar{F}_0^{-1}(\bar{F}_0(t)\bar{F}_0(x_i)))} \frac{\partial  (\bar{F}_0^{-1}(\bar{F}_0(t)\bar{F}_0(x_i)))}{\partial t} = \theta r_0(t)
	\end{align*}
	which gives a condition in terms of the hazard gradient to satisfy the functional equation \eqref{kundu_fe_gen} as
	\begin{equation}\label{kundu_hr_grad}
		\sum_{i=1}^{2}r_i(\bar{F}_0^{-1}(\bar{F}_0(t)\bar{F}_0(x_1)), \bar{F}_0^{-1}(\bar{F}_0(t)\bar{F}_0(x_2))) \frac{\partial (\bar{F}_0^{-1}(\bar{F}_0(t)\bar{F}_0(x_i)))}{\partial t} = \theta r_0(t),
	\end{equation}
	where $r_0(\cdot)$ is the baseline hazard function. On integrating \eqref{kundu_hr_grad}, we can retrieve the functional equation \eqref{kundu_fe_gen}.

	\begin{thm}
		A bivariate random vector $X=(X_1,X_2)$ satisfy \eqref{kundu_fe_gen} if and only if
		\begin{align}\label{kundu_hr_grad_marg}
			r_i(x_1,x_2)=
			\begin{cases}
				{r_{i} \bigg( \bar{F}_0^{-1} \Big( \frac{\bar{F}_0(x_i)}{\bar{F}_0(x_{3-i})} \Big) \bigg) \frac{\partial}{\partial x_{i}} \bigg( \bar{F}_0^{-1} \Big( \frac{\bar{F}_0(x_i)}{\bar{F}_0(x_{3-i})} \Big) \bigg) } & ;~ x_i > x_{3-i} \\
				{-r_{3-i} \bigg( \bar{F}_0^{-1} \Big( \frac{\bar{F}_0(x_{3-i})}{\bar{F}_0(x_{i})} \Big) \bigg) \bigg \vert \frac{\partial}{\partial x_{i}} \bigg( \bar{F}_0^{-1} \Big( \frac{\bar{F}_0(x_{3-i})}{\bar{F}_0(x_{i})} \Big) \bigg) \bigg \vert +  \theta r_0(x_i) } & ;~ x_i < x_{3-i},
			\end{cases}
		\end{align}	
		where $r_i(x)$ is the marginal hazard rate of $X_i;~i=1,2$.
	\end{thm}
	
	\begin{proof}
		Suppose $(X_1,X_2)$ satisfy,
		\begin{align*}
			\bar{F}(\bar{F}_0^{-1}(\bar{F}_0(t)\bar{F}_0(x_1)), \bar{F}_0^{-1}(\bar{F}_0(t)\bar{F}_0(x_2)))=\bar{F}(x_1,x_2) [\bar{F}_0(t)]^{\theta} \\
			\implies \bar{F}(x_1,x_2) = \bar{F} \bigg( \bar{F}_0^{-1} \Big( \frac{\bar{F}_0(x_1)}{\bar{F}_0(t)} \Big), \bar{F}_0^{-1} \Big( \frac{\bar{F}_0(x_2)}{\bar{F}_0(t)} \Big) \bigg) [\bar{F}_0(t)]^{\theta}
		\end{align*}
		For $t=\text{min}\{x_1,x_2\}$ we have for $x_1>x_2$,
		\begin{align*}
			\bar{F}(x_1,x_2)&=\bar{F} \bigg( \bar{F}_0^{-1} \Big( \frac{\bar{F}_0(x_1)}{\bar{F}_0(x_2)} \Big), 0 \bigg) [\bar{F}_0(x_2)]^{\theta}\\
			&=\bar{F}_1 \bigg( \bar{F}_0^{-1} \Big( \frac{\bar{F}_0(x_1)}{\bar{F}_0(x_2)} \Big) \bigg) [\bar{F}_0(x_2)]^{\theta}.
		\end{align*}
		Taking negative of the logarithms on both sides and differentiating w.r.t $x_i;~i=1,2$ we get,
		$$r_1(x_1,x_2)=r_1 \bigg( \bar{F}_0^{-1} \Big( \frac{\bar{F}_0(x_1)}{\bar{F}_0(x_2)} \Big) \bigg) \frac{\partial}{\partial x_1} \bigg( \bar{F}_0^{-1} \Big( \frac{\bar{F}_0(x_1)}{\bar{F}_0(x_2)} \Big) \bigg) $$
		and $$r_2(x_1,x_2)=-r_1 \bigg( \bar{F}_0^{-1} \Big( \frac{\bar{F}_0(x_1)}{\bar{F}_0(x_2)} \Big) \bigg) \bigg \vert \frac{\partial}{\partial x_2} \bigg( \bar{F}_0^{-1} \Big( \frac{\bar{F}_0(x_1)}{\bar{F}_0(x_2)} \Big) \bigg) \bigg \vert + \theta r_0(x_2).$$
		In a similar manner we can prove the result for $x_1 < x_2$ to obtain \eqref{kundu_hr_grad_marg}. The converse is straight forward through the representation 
		\begin{equation*}
			\bar{F}(x_1,x_2)=exp\Bigg [-\int_{0}^{x_1}r_1(u,0)du-\int_{0}^{x_2}r_2(x_1,u)du\Bigg ]
		\end{equation*}
		or
		\begin{equation*}
			\bar{F}(x_1,x_2)=exp\Bigg [-\int_{0}^{x_1}r_1(u,x_2)du-\int_{0}^{x_2}r_2(0,u)du\Bigg ]
		\end{equation*}
		(\citet{johnson1975vector}).
	\end{proof}
	If we closely investigate equation \eqref{kundu_hr_grad_marg}, we see that the marginal hazard function $r_i(x)$ should satisfy the condition
	\begin{equation}\label{kundu_marg_hr_bound}
		r_{i} \Big( \bar{F}_0^{-1} \Big( \frac{\bar{F}_0(x_{i})}{\bar{F}_0(x_{3-i})} \Big) \Big) \bigg \vert \frac{\partial}{\partial x_{3-i}} \Big( \bar{F}_0^{-1} \Big( \frac{\bar{F}_0(x_{i})}{\bar{F}_0(x_{3-i})} \Big) \Big) \bigg \vert \leq  \theta r_0(x_{3-i});~\text{for $i=1,2$}
	\end{equation}
	since $r_i(x_1,x_2) \geq 0$. Based on these observations necessary and sufficient conditions for differentiable $r_i(x)$ to qualify as marginal failure rates for $\bar{F}(x_1,x_2)$ to be a bivariate survival function is discussed in the next theorem. 
	
	\begin{thm}\label{kundu_marg_hr_suffcond}
		The necessary and sufficient conditions for differentiable functions $r_{i}(x);~i=1,2$ to qualify as marginal failure rates of $\bar{F}(x_1,x_2)$ in \eqref{kundu_surv_gen} is that for $x_i > x_{3-i};~i=1,2$,
		\begin{enumerate}[label=(\roman*)]
			\item $0 \leq r_{i} \Big( \bar{F}_0^{-1} \big( \frac{\bar{F}_0(x_{i})}{\bar{F}_0(x_{3-i})} \big) \Big) \bigg \vert \frac{\partial \bar{F}_0^{-1} \big( \frac{\bar{F}_0(x_{i})}{\bar{F}_0(x_{3-i})} \big)}{\partial x_{3-i}} \bigg \vert \leq  \theta r_0(x_{3-i})$ \label{kundu_hr_cond1.1}
			\item $\int_{0}^{\infty}r_i(x)dx=\infty$ \label{kundu_hr_cond1.2}
			\item $r_i\Big( \bar{F}_0^{-1} \big( \frac{\bar{F}_0(x_{i})}{\bar{F}_0(x_{3-i})} \big) \Big) \frac{\partial \bar{F}_0^{-1} \big( \frac{\bar{F}_0(x_{i})}{\bar{F}_0(x_{3-i})} \big)}{\partial x_i} \Big[\theta r_0(x_{3-i}) + r_i\Big( \bar{F}_0^{-1} \big( \frac{\bar{F}_0(x_{i})}{\bar{F}_0(x_{3-i})} \big) \Big) \frac{\partial \bar{F}_0^{-1} \big( \frac{\bar{F}_0(x_{i})}{\bar{F}_0(x_{3-i})} \big)}{\partial x_{3-i}} \Big]-$ \\ $\frac{\partial}{\partial x_{3-i}} [r_i\Big( \bar{F}_0^{-1} \big( \frac{\bar{F}_0(x_{i})}{\bar{F}_0(x_{3-i})} \big) \Big) \frac{\partial \bar{F}_0^{-1} \big( \frac{\bar{F}_0(x_{i})}{\bar{F}_0(x_{3-i})} \big)}{\partial x_i}] \geq 0$ \label{kundu_hr_cond1.3}
			\item $\theta \leq v_1(\Theta) + v_2(\Theta) \leq  2 \theta$ \label{kundu_hr_cond1.4}
		\end{enumerate}
		where $v_i(\Theta)=\Big[ r_i\Big( \bar{F}_0^{-1} \big( \frac{\bar{F}_0(x_{i})}{\bar{F}_0(x_{3-i})} \big) \Big) exp\{ - \int_{0}^{\bar{F}_0^{-1} \big( \frac{\bar{F}_0(x_{i})}{\bar{F}_0(x_{3-i})} \big)} r_i(u) du \} \Big \vert \frac{\partial \bar{F}_0^{-1} \big( \frac{\bar{F}_0(x_{i})}{\bar{F}_0(x_{3-i})} \big)}{\partial (- \ln \bar{F}_0(x_{3-i}))} \Big \vert \Big ] _{\vert x_i=x_{3-i}};~i=1,2$. The bivariate function $\bar{F}(x_1,x_2)$ given by equation \eqref{kundu_surv_gen} is a survival function satisfying the functional equation in \eqref{kundu_fe_gen} where
		\begin{equation}\label{kundu_marg_hr_suffcond_surv}
			\bar{F}_i(x)=exp \big( -\int_{0}^{x} r_i(u) du \big);~x \geq 0,~i=1,2.
		\end{equation}
	\end{thm}
	
	\begin{proof}
		Condition \ref{kundu_hr_cond1.1} follows from equation \eqref{kundu_marg_hr_bound} and condition \ref{kundu_hr_cond1.2} is necessary for $r_i;~i=1,2$ to be a hazard rate function. Conditions \ref{kundu_hr_cond1.3} and \ref{kundu_hr_cond1.4} implies conditions \ref{kundu_marg_densty1} and \ref{kundu_marg_densty2} of Theorem \ref{kundu_mo_theorem}.
		\par Conversely the hazard gradient associated with $\bar{F}(x_1,x_2)$ is given by \eqref{kundu_hr_grad_marg} and the condition \ref{kundu_hr_cond1.1} follows from the non-negativity of $r_i(x_1,x_2);~i=1,2$. Since $r_i(\cdot);~i=1,2$ are univariate failure rates, condition \ref{kundu_hr_cond1.2} follows. Condition \ref{kundu_hr_cond1.3} follows from the non-negativity of the density function associated with $\bar{F}(x_1,x_2)$ and condition \ref{kundu_hr_cond1.4} follows from Theorem \ref{kundu_mo_theorem}.
	\end{proof}

	
	\textbf{Counter Example \ref{kundu_counter_eg} (contd.)} LFR-exponential: The marginal failure rates are given by
	\begin{equation*}
		r_i(x_i)=2\alpha x_i + 1;~i=1,2.
	\end{equation*}
	Then the condition \ref{kundu_hr_cond1.1} in Theorem \ref{kundu_marg_hr_suffcond} becomes $0 \leq 2 \alpha
	(x_i - x_{3-i}) + 1 \leq \theta;~x_i > x_{3-i}$ which is not satisfied for the choice of $\alpha=1.5,~\theta=3,~x_1=5$ and $x_2=3$. Hence, we cannot construct a bivariate survival function $\bar{F}(x_1,x_2)$ with the marginals as linear failure rate and baseline distribution as exponential.
	
	\begin{eg}
		\normalfont BPHC-Weibull: For a proportional hazard rate, $$r_i(x_i)=(\theta_i + \theta_3) r_0(x_i);~x_i > 0,~i=1,2$$ with $\theta_i > 0;~i=1,2,3$ and $r_0(x)=\alpha x^{\alpha - 1};~x > 0,~\alpha > 0$, it is easily seen that conditions of Theorem \ref{kundu_marg_hr_suffcond} is satisfied. For
		\begin{enumerate}[label=(\roman*)]
			\item $0 \leq \theta_i + \theta_3 \leq \theta_1 + \theta_2 + \theta_3$
			\item $\int_{0}^{\infty} \alpha (\theta_i + \theta_3)x_{i}^{\alpha-1} dx_i = \infty$
			\item $\theta_i + \theta_3 \geq 0$
			\item $\theta_1 + \theta_2 + \theta_3 \leq \theta_1 + \theta_2 + 2\theta_3 \leq 2(\theta_1 + \theta_2 + \theta_3)$.
		\end{enumerate}
		for $x_i > x_{3-i};~i=1,2$. From \eqref{kundu_marg_hr_suffcond_surv}, $\bar{F}_i(x_i)=[\bar{F}_0(x_i)]^{\theta_i + \theta_3};~i=1,2$, where $\bar{F}_0(x)=e^{-x^{\alpha}};~x > 0,~\alpha > 0$. Substituting in \eqref{kundu_surv_gen}, we get
		\begin{align*}\label{kundu_surv_weibull}
			\bar{F}(x_1,x_2)&=
			\begin{cases}
				e^{-(\theta_1 + \theta_3)x_1^{\alpha} - \theta_2 x_2^{\alpha}} &~;~x_{1} \geq x_{2}\\
				e^{-\theta_1 x_1^{\alpha} - (\theta_2 + \theta_3) x_2^{\alpha}} &~;~x_{1} \leq x_{2}.
			\end{cases}
		\end{align*}
	\end{eg}
	
	Hence a class of distributions which include the bivariate semi-parametric singular family of distributions given by \citet{kundu2022bivariate} has been proposed as a general solution to the functional equation in \eqref{kundu_fe_gen}. The estimation of popular distributions belonging to this class have been already discussed in the literature (see \citet{proschan1976estimating}, \citet{pena1990bayes}, \citet{kundu2022bivariate}).
	
	\subsection*{\textbf{Conflict of Interest Statement}}
	The authors have no conflicts of interest to declare.
	
	\subsection*{\textbf{Acknowledgments}}
	This work was financially supported by the Council of Scientific \& Industrial Research (CSIR), Government of India through the Senior Research Fellowship scheme vide No 09/239(0551) /2019-EMR-I.
	
	\bibliography{ref}
	
\end{document}